\documentclass[12pt]{article}
\usepackage{MnSymbol}
\usepackage{amsthm}

\DeclareMathAlphabet\mathbb{U}{msb}{m}{n}
\DeclareMathAlphabet\Bbb{U}{msb}{m}{n}

\usepackage{vmargin}
\usepackage{wasysym}
\usepackage[active]{srcltx}
\usepackage[stable]{footmisc}
\usepackage{caption}
\usepackage{pdfpages}
\usepackage[matrix, arrow,all,cmtip,color]{xy}
\usepackage{url}
\usepackage{epigraph}

\newcommand{\bi}{\begin{itemize}}
\newcommand{\ei}{\end{itemize}}
\newcommand{\bd}{\begin{description}}
\newcommand{\ed}{\end{description}}

\theoremstyle{plain}
\newtheorem{conjecture}{Conjecture}

\newtheorem{claim}{Claim}

\theoremstyle{definition}

\newtheorem{remark}{Remark}
\newtheorem{problem}{Question}

\def\bqqq{\begin{quote}}
\def\eqqq{\end{quote}}

\def\lra{\longrightarrow}
\def\ra{\rightarrow}
\def\rtt{\,\rightthreetimes\,}
\def\xra{\xrightarrow}

\def\ZZ{\Bbb Z}

\def\rrt#1#2#3#4#5#6{\xymatrix{ {#1} \ar[r]^{} \ar@{->}[d]_{#2} & {#4} \ar[d]^{#5} \\ {#3}  \ar[r] \ar@{-->}[ur]^{}& {#6} }}

\def\union{\cup}

\def\RR{\Bbb R}

\def\card{\,{\mathrm{card}\,}}

\begin{document}
\date{9 May 2017}
\title{
The unreasonable power of the lifting property in elementary mathematics
}
\author{
misha gavrilovich\thanks{A draft; comments welcome. {\tt mi\!\!\!ishap\!\!\!p@sd\!\!\!df.org. current version at http://mishap.sdf.org/expressive-power-of-the-lifting-property.pdf.} Minor update 7.17.
}\\
in memoriam: evgenii shurygin
}
\maketitle
\setlength{\epigraphwidth}{0.8\textwidth}
\epigraph{
   instances of human and animal behavior
[...]  miraculously complicated, 
[...] they have
  little, if any, pragmatic (survival/reproduction) value.
[...]   they are due to internal constraints on
   possible architectures of unknown to us functional "mental structures".
}{
   Gromov, Ergobrain
}

\begin{abstract}

We illustrate the generative power of the lifting property (orthogonality 
of morphisms in a category) 
as a means of
defining natural elementary mathematical concepts by giving a number of examples in
various categories, in particular showing that many standard elementary notions
of abstract topology can be defined by applying the lifting property to simple
morphisms of finite topological spaces. 
Examples in topology include the notions of: compact, discrete, connected, and
totally disconnected spaces, dense image, induced topology, and separation axioms.
Examples in algebra include: finite groups being 
nilpotent, solvable, torsion-free, $p$-groups, and prime-to-$p$ groups;
injective and projective
 modules; injective, surjective,
and split homomorphisms.

We include some speculations on the
wider significance of this.
\end{abstract}

\section{
Introduction.
}

The purpose of this short note is to draw attention to the following
observation which we find rather curious:
\begin{quote}
    a number of elementary properties from a first-year course
    can be defined category-theoretically by
    repeated application of a standard category theory trick,
    the Quillen lifting property,
    starting from a class of explicitly given morphisms,
    often consisting of a single (counter)example
\end{quote}

In particular,  several elementary notions of topology
have Kolmogorov complexity of several bytes
in a natural category-theoretic formalism (explained below),
e.g. \begin{quote}
     compactness is   "$((\{a\}\longrightarrow \{a{\small\searrow}b\})^r_{<5})^{lr}$",  \\
     connectedness is "$(\{a,b\}\longrightarrow \{a=b\})^l$",       \\
     dense image is     "$(\{b\}\longrightarrow \{a{\small\searrow}b\})^l$".      \\
     \end{quote}

We suggest it appears worthwhile to try to develop a formalism,
or rather a very short program (kilobytes) based on such a formalism,
which supports reasoning in elementary topology.

These observations arose in an attempt to understand ideas of
Misha Gromov [Memorandum Ergo] about ergologic/ergostructure/ergosystems.
Oversimplifying, ergologic is a kind of reasoning which helps to
understand how to
generate proper concepts, ask interesting questions, and, more generally,
produce interesting rather than useful or correct behaviour.
He conjectures there is a related class of
mathematical, essentially combinatorial, structures,
called {\em ergostructures} or {\em ergosystems},
and that this concept might eventually
help to understand complex biological behaviour including learning
and create  mathematically interesting models of these processes.

We hope our observations may eventually help to uncover
an essentially combinatorial reasoning behind elementary topology,
and thereby suggest an example of an ergostructure.

{\em Related works.} 
This paper continues work started in [DMG], a rather leisurely 
introduction to some of the ideas presented here.  
Draft [Gavrilovich, Elementary Topology] shows how to view several
topology notions and arguments in [Bourbaki, General Topology] 
as diagram chasing calculations with finite categories.
Draft [Gavrilovich, Tame Topology] is more speculative but less verbose;
it has several more examples dealing with compactness, 
in particular it shows that a number of consequences of compactness 
can be expressed as a change of order of quantifiers in a formula. 
Notably, these drafts show
how to "read off" a simplicial topological space 
from the definition of a uniform space, see also Remarks~\ref{sTop:1}~and~\ref{simp:obj:unifrom}.

{\em Structure of the paper.} A mathematically inclined reader might want to read
only the first two sections with miscellaneous examples of lifting properties
and a combinatorial notation for elementary properties of topological spaces.
A logician might want to read in the third section our suggestions towards 
a theorem prover/proof system for elementary topology based on diagram
chasing. Appendix~A states separation axioms in terms of lifting properties
and finite topological spaces. Appendix~B reproduces some references we use. 

Finally, the last section attempts to explain our motivation and says
a few words about the concept of ergostructure by Misha Gromov.

We would also like to draw attention to Conjecture~1 (a charaterisation of the class of proper maps) 
and Question~2 asking for a characterisation of the circle and the interval.

\section{
The lifting property: the key observation
}

For a property \ensuremath{C} of arrows (morphisms) in a category, define

$$ C^l := \{ \ensuremath{f} :\text{ for each }g \in C\ \ensuremath{f} \,\rightthreetimes\,  \ensuremath{g} \} $$
$$ C^r := \{ \ensuremath{g} :\text{ for each }f \in C\ \ensuremath{f} \,\rightthreetimes\,  \ensuremath{g} \} $$
$$ C^{lr}:=(C^l)^r, ... $$

here $f \,\rightthreetimes\,  g$ reads " $f$ has the left lifting property wrt $g$ ",
" $f$ is (left) orthogonal to $g$ ",
i.e.  for  $f:A\longrightarrow B$, $g:X\longrightarrow Y$,
$f \,\rightthreetimes\, g$ iff for each $i:A\longrightarrow X$, $j:B\longrightarrow Y$ such that $ig=fj$ ("the square commutes"),
there is $j':B\longrightarrow X$ such that $fj'=i$ and $j'g=j$ ("there is a diagonal
making the diagram commute").

The following observation is enough to reconstruct all the examples in this
paper, with a bit of search and computation.
\begin{quote}
 {\bf Observation.}\\ 
  A number of elementary properties can be obtained by repeatedly passing
  to the left or right orthogonal $C^l, C^r, C^{lr}, C^{ll}, C^{rl}, C^{rr},...$
  starting from a simple class of morphisms, often
  a single (counter)example to the property you define.
\end{quote}

A useful intuition is to think that the property of left-lifting against a
class \ensuremath{C} is a kind of negation of the property of being in \ensuremath{C}, and that
right-lifting is another kind of negation.  Hence the classes obtained from C
by taking orthogonals an odd number of times, such as $C^l, C^r, C^{lrl}, C^{lll}$
etc., represent various kinds of negation of \ensuremath{C}, so $C^l, C^r, C^{lrl}, C^{lll}$ each
consists of morphisms which are far from having property C.

Taking the orthogonal of a class \ensuremath{C} is a simple way to define a class of
morphisms excluding non-isomorphisms from \ensuremath{C}, in a way which is useful in a
diagram chasing computation.

The class $C^l$ is always closed under retracts, pullbacks, (small) products
(whenever they exist in the category) and composition of morphisms, and
contains all isomorphisms of C. Meanwhile, $C^r$ is closed under retracts,
pushouts, (small) coproducts and transfinite composition (filtered colimits) of
morphisms (whenever they exist in the category), and also contains all
isomorphisms.    

For example, the notion of isomorphism can be obtained starting from the class
of all morphisms, or any single example of an isomorphism:
$$
(Isomorphisms) = (all\ morphisms)^l = (all\ morphisms)^r = (h)^{lr} = (h)^{rl}
$$
where \ensuremath{h } is an arbitrary isomorphism.

\def\rrt#1#2#3#4#5#6{\xymatrix{ {#1} \ar[r]^{} \ar@{->}[d]_{#2} & {#4} \ar[d]^{#5} \\ {#3}  \ar[r] \ar@{-->}[ur]^{}& {#6} }}
\begin{figure}
\begin{center}
\large
$ (a)\ \xymatrix{ A \ar[r]^{i} \ar@{->}[d]_f & X \ar[d]^g \\ B \ar[r]_-{j} \ar@{-->}[ur]^{{\tilde j}}& Y }$
$(b)\  \rrt  {\{\}}  {} {\{\bullet\}}  X {\therefore(surj)} Y $
$(c)\  \rrt {\{\bullet,\bullet\}} {} {\{\bullet\}}  X {\therefore(inj)} Y $
\end{center}
\caption{\label{fig1}\normalsize
Lifting properties. 
 (a) The definition of a lifting property $f\rtt g$. 
 (b) $X\lra Y$ is surjective 
 (c) $X\lra Y$ is injective 
}
\end{figure}

{\bf 
Example.} \begin{quote}
  Take $C=\{ \emptyset\longrightarrow \{*\} \}$ in Sets and Top. 
  Let us show that $C^l$ is the class of surjections, $C^{rr}$ is the class of subsets, 
  $C^l$ consists of maps $f:A\longrightarrow B$ such that either  $A=B=\emptyset$ 
  or $A\neq \emptyset$, \ensuremath{B} arbitrary. Further, 
  in Sets, $C^{rl}$ is the class of injections, and in Top, $C^{rl}$ 
  is the class of maps of form $A\longrightarrow A\cup D$, $D$ is discrete.

   $A\longrightarrow \ensuremath{B} \,\rightthreetimes\,  \emptyset\longrightarrow \{*\}$ iff $A=B=\emptyset$ or $A\neq \emptyset$, \ensuremath{B} arbitrary. Indeed, 
  if $A\neq \emptyset$, there is no map $i:A\longrightarrow X=\emptyset$ and the lifting property 
  holds vacuously; if $A=\emptyset\neq B$, there exist unique maps $i:A=\emptyset\longrightarrow X=\emptyset$, 
  $j:B\longrightarrow Y=\{*\}$, but no map $j':B\longrightarrow X=\emptyset$ as $B\neq \emptyset$ by assumption. 

  $\emptyset\longrightarrow \{*\} \,\rightthreetimes\,  g$ iff $g:X\longrightarrow Y$ is surjective; indeed, 
  the map $j:B=\{*\}\longrightarrow Y$ picks a point in $Y$ and $j':B=\{*\}\longrightarrow X$
  picks its preimage as $j'g=j$; the other condition $fj'=i:\emptyset\longrightarrow X$ 
  holds trivially. Thus $(\emptyset\longrightarrow \{*\})^r$ is the class of surjections.  
  
  In Sets, $(\emptyset\longrightarrow \{*\})^{rl}$ is the class of injections, i.e.
   $f \,\rightthreetimes\,  g$ for each surjection \ensuremath{g} iff \ensuremath{f} is injective; indeed,
  for such \ensuremath{f} and \ensuremath{g} the following is well-defined:
  set $j'(b)=i(f^{ -1}(b))$ for \ensuremath{b} in $Im f$, and $j'(b)=g^{ -1}(j(b))$ 
  otherwise; for injective \ensuremath{f} $j'(b)$ does not depend on the choice 
  of a preimage of b, and for \ensuremath{g} surjective a preimage always exists. 

  In Top, $(\emptyset\longrightarrow \{0\})^{rl}$ is the class of maps of form $A\longrightarrow A\cup D$, $D$ is discrete;
  given a map $A\longrightarrow B$, consider $A\longrightarrow \ensuremath{B} \,\rightthreetimes\,  ImA\cup (B\setminus A)\longrightarrow B$ where $ImA\cup D\longrightarrow B$
  denotes the disjoint union of the image of A in \ensuremath{B} with induced topology,
  and $B\setminus A$ equipped with the discrete topology. 

  In both Sets and Top,  $(\emptyset\longrightarrow \{*\})^{rr}$ is the class of subsets, i.e. 
  injective maps $A\hookrightarrow B$ where the topology on $A$ is induced from $B$. 
\end{quote}

Toying with the observation leads to the examples in the claim below
which is trivial to verify, an exercise in deciphering the notation
in all cases but (vii) proper.

\begin{claim}\begin{itemize}
\item[ (i)] $(\emptyset\longrightarrow \{*\})^r$, $(0\longrightarrow R)^r$, and $\{0\longrightarrow \ZZ\}^r$ are the classes of surjections in
      in tha categories of Sets, R-modules, and Groups, resp.,
      (where $\{*\}$ is the one-element set, and in the category of (not necessarily abelian) groups, $0$ denotes the trivial group)
\item[ (ii)] $(\{\star,\bullet\}\longrightarrow \{*\})^l=(\{\star,\bullet\}\longrightarrow \{*\})^r$, $(R\longrightarrow 0)^r$, $\{\ZZ\longrightarrow 0\}^r$ are  the classes of
      injections in the categories of Sets, R-modules, and Groups, resp
\item[ (iii)] in the category of R-modules, \\
                 \ \ a module \ensuremath{P} is projective iff $0\longrightarrow P$ is in $(0\longrightarrow R)^{rl}$\\
                 \ \  a module I is injective iff $I\longrightarrow 0$ is in $(R\longrightarrow 0)^{rr}$
\item[ (iv)] in the category of Groups, \begin{itemize}
\item[]        a finite group \ensuremath{H} is nilpotent iff $H\longrightarrow H\times H$ is in $\{\, 0\longrightarrow \ensuremath{G} : G\text{ arbitrary} \}^{lr}$ 
\item[]        a finite group \ensuremath{H} is solvable iff $0\longrightarrow H$ is in $\{\, 0\longrightarrow A : A\text{ abelian }\}^{lr}= \{\, [G,G]\longrightarrow \ensuremath{G} : G\text{ arbitrary }\}^{lr}$
\item[]        a finite group \ensuremath{H} is of order prime to $p$ iff $H\longrightarrow 0$ is in $\{\ZZ/p\ZZ\longrightarrow 0\}^r$
\item[]        a finite group \ensuremath{H} is a p-group iff $H\longrightarrow 0$ is in $\{\ZZ/p\ZZ\longrightarrow 0\}^{rr}$
\item[]        a group \ensuremath{H} is torsion-free iff $0\longrightarrow H$ is in $\{ n\ZZ\longrightarrow \ZZ: n>0 \}^r$
\item[]        a group \ensuremath{F} is free iff $0\longrightarrow F$ is in  $\{0\longrightarrow \ZZ\}^{rl}$
\item[]        a homomorphism $f$ is split iff $f \in  \{\, 0\longrightarrow \ensuremath{G} : G\text{ arbitrary} \}^r$
\end{itemize}
\item[ (v)] in the category of metric spaces and uniformly continuous maps,\\
        a metric space \ensuremath{X} is complete iff $\{1/n\}_n\longrightarrow \{1/n\}_n\cup \{0\} \,\rightthreetimes\,  X\longrightarrow \{0\}$
           where the metric on $\{1/n\}_n$ and $\{1/n\}_n\cup \{0\}$ is induced from the
                 real line\\
        a subset $A \subset  X$ is closed iff  $\{1/n\}_n\longrightarrow \{1/n\}_n\cup \{0\} \,\rightthreetimes\,  A\longrightarrow X$
\item[ (vi)]    in the category of topological spaces,\\
     for a connected topological space X, each function on \ensuremath{X} is bounded
     iff 
         $$ \emptyset\longrightarrow \ensuremath{X} \,\rightthreetimes\,  \cup_n (-n,n) \longrightarrow  \RR$$ 
\item[ (vii)] in the category of topological spaces (see notation defined below),
\begin{itemize}

\item[] a Hausdorff space $K$ is compact iff $K\longrightarrow \{*\}$ is in  $((\{a\}\longrightarrow \{a{\small\searrow}b\})^r_{<5})^{lr}$
\item[] a  Hausdorff space $K$ is compact iff $K\longrightarrow \{*\}$ is in  $$
     \{\, \{a\leftrightarrow b\}\longrightarrow \{a=b\},\, \{a{\small\searrow}b\}\longrightarrow \{a=b\},\,
     \{b\}\longrightarrow \{a{\small\searrow}b\},\,\{a{\small\swarrow}o{\small\searrow}b\}\longrightarrow \{a=o=b\}\,\,\}^{lr}$$

\item[] a space $D$ is discrete iff $ \emptyset \longrightarrow  D$ is in $   (\emptyset\longrightarrow \{*\})^{rl}      $  

\item[] a space $D$ is antidiscrete iff $ \ensuremath{D} \longrightarrow  \{*\} $ is in 
$(\{a,b\}\longrightarrow \{a=b\})^{rr}= (\{a\leftrightarrow b\}\longrightarrow \{a=b\})^{lr} $ 

\item[] a space $K$ is connected or empty iff $K\longrightarrow \{*\}$ is in  $(\{a,b\}\longrightarrow \{a=b\})^l $
\item[] a space $K$ is totally disconnected and non-empty iff $K\longrightarrow \{*\}$ is in  $(\{a,b\}\longrightarrow \{a=b\})^{lr} $

\item[] a space $K$ is connected and non-empty  
 iff 
 for some arrow $\{*\}\longrightarrow K$\\
$\text{ \ \ \ \ \     } \{*\}\longrightarrow K$ is in 
            $   (\emptyset\longrightarrow \{*\})^{rll} = (\{a\}\longrightarrow \{a,b\})^l$ 
       
\item[] a space $K$ is non-empty iff $K\longrightarrow \{*\}$ is in $   (\emptyset\longrightarrow \{*\})^l$ 
\item[] a space $K$ is empty iff $K \longrightarrow \{*\}$ is in $   (\emptyset\longrightarrow \{*\})^{ll}$ 
\item[] a space $K$ is $T_0$ iff $K  \longrightarrow \{*\}$ is in $   (\{a\leftrightarrow b\}\longrightarrow \{a=b\})^r$  
\item[]  a space $K$ is $T_1$ iff $K  \longrightarrow \{*\}$ is in $   (\{a{\small\searrow}b\}\longrightarrow \{a=b\})^r$ 
\item[] a space $X$ is Hausdorff iff for each injective map $\{x,y\} \hookrightarrow  X$ 
it holds $\{x,y\} \hookrightarrow  \ensuremath{X} \,\rightthreetimes\,  \{ \ensuremath{x} {\small\searrow} \ensuremath{o} {\small\swarrow} \ensuremath{y} \} \longrightarrow  \{ x=o=y \}$

\item[] a non-empty space $X$ is regular (T3) iff for each arrow $    \{x\} \longrightarrow  X$ it holds 
    $    \{x\} \longrightarrow  \ensuremath{X} \,\rightthreetimes\,  \{x{\small\searrow}X{\small\swarrow}U{\small\searrow}F\} \longrightarrow  \{x=X=U{\small\searrow}F\}$
\item[] a space $X$ is normal (T4) iff $\emptyset \longrightarrow \ensuremath{X} \,\rightthreetimes\,   \{a{\small\swarrow}U{\small\searrow}x{\small\swarrow}V{\small\searrow}b\}\longrightarrow \{a{\small\swarrow}U=x=V{\small\searrow}b\}$ 

\item[] a space $X$ is completely normal iff $\emptyset\longrightarrow \ensuremath{X} \,\rightthreetimes\,  [0,1]\longrightarrow \{0{\small\swarrow}x{\small\searrow}1\}$ 
 where the map $[0,1]\longrightarrow \{0{\small\swarrow}x{\small\searrow}1\}$ sends $0$ to $0$, $1$ to $1$, and the rest $(0,1)$ to $x$ 
\item[] a space $X$ is path-connected iff $\{0,1\} \longrightarrow  [0,1] \,\rightthreetimes\,  \ensuremath{X} \longrightarrow  \{*\}$ 
\item[] a space $X$ is path-connected iff for each Hausdorff compact space $K$ and each injective map $\{x,y\} \hookrightarrow  K$ it holds
   $\{x,y\} \hookrightarrow  \ensuremath{K} \,\rightthreetimes\,  \ensuremath{X} \longrightarrow  \{*\}$

\item[]       $(\emptyset\longrightarrow \{*\})^r$   is the class of surjections
\item[]       $(\emptyset\longrightarrow \{*\})^{rr}$ is the class of subsets, i.e. injective maps $A\hookrightarrow B$ where the topology on $A$ is induced from $B$
\item[]       $(\emptyset\longrightarrow \{*\})^{lll}$ is the class of maps $A\longrightarrow B$ which split

\item[]       $(\{b\}\longrightarrow \{a{\small\searrow}b\})^l$ is the class of maps with dense image
\item[]       $(\{b\}\longrightarrow \{a{\small\searrow}b\})^{lr}$ is the class of closed subsets $A \subset  X$, $A$ a closed subset of $X$
\item[]       $((\{a\}\longrightarrow \{a{\small\searrow}b\})^r_{<5})^{lr}$ is roughly the class of proper maps
       (see below).
\end{itemize}
\end{itemize}
\end{claim}
{\bf
Proof.
} In (iv), we use that a finite group $H$ is nilpotent iff the diagonal $\{
(h,h) : \ensuremath{h } \in \ensuremath{H} \}$ is subnormal in $H\times H$. (vii) is discussed below. QED.

Appendix~A shows that the usual formulations of the separation axioms are in fact lifting properties.

\begin{problem} 
Find more examples of meaningful lifting properties in various categories.
        Play with natural classes of morphisms to see whether their iterated orthogonals         
        are meaningful.
\end{problem}

\begin{remark} Most of the definitions above are in a form useful in a diagram chasing computation.
Let us comment on this.

In category theory it is usual to view an object $X$ as either 
the arrow $\perp\longrightarrow X$ or $X\lra\top$ from the initial object $\perp$ or to the terminal object $\top$.
However, in Claim~1(vii) we sometimes view (the properties of) a space as (properties of a non-unique!) arrow $\{*\} \longrightarrow X$ or
$\{a,b\} \hookrightarrow X$. Our purpose is to give definitions useful in a diagram chasing computation,
and these definitions can be used in this way. 

A group $G$ has cohomological dimension \ensuremath{1} iff each surjections $G'\longrightarrow G$ splits. Item (iv) views this 
as the following diagram chasing rule: 
 in (the valid diagram corresponding to)  $0\longrightarrow \ZZ \,\rightthreetimes\,  G'\longrightarrow G$, it is permissible to replace $\ZZ$ by 
 an arbitrary group $A$ (thereby obtaining a valid diagram corresponding to $0\longrightarrow A \,\rightthreetimes\,  G'\longrightarrow G$). 

\end{remark}

\begin{remark}\label{sTop:1}
Claim (v) shows that a computer-generated proof in (Ganesalingam, Gowers; Problem 2)
        of the claim that completeness is inherited by closed subsets of metric spaces, i.e. 
        a closed subspace of a complete metric space is necessarily complete,
        translates to  two applications of a diagram chasing rule
        corresponding to the lifting property. 

In fact, in any category for an arbitrary class $C$ of morphisms it holds
 $ X\longrightarrow \{*\} \,\in\, C^l$ and $f: A \hookrightarrow  \ensuremath{X} \,\in\, C^l $ implies 
$A \longrightarrow  \{*\}  \,\in\, C^l $ whenever $f$ is a monomorphism and where $\{*\}$ denotes the terminal object. 
\end{remark}

Note that, 
the definition involved infinite objects 
which are infinite sequences
$\{1/n\}_n\longrightarrow \{1/n\}_n\cup \{0\} $. This can be probably be avoided if instead we 
consider {\em uniform spaces} as {\em simplicial objects in the category 
of topological spaces} (see Remark~\ref{simp:obj:unifrom} for details)
and interpret the lifting property in a subcategory of the category of simplicial objects
in the category of topological spaces using that a Cauchy filter is (almost) a map
from a ``constant'' simplicial object. 

This would give a formal analogy corresponding to
the informal analogy between {\em compact} topological spaces and 
{\em complete} metric (uniform) spaces.  

\begin{remark}
Claim demonstrates that a number of elementary notions can be concisely
        expressed in terms of a simple diagram chasing rule. However,
        it appears there is no (well-known) logic or proof system
        based on diagram chasing in a category. We make some suggestions towards
        such a proof system in the last two sections.
\end{remark}

\section{
Elementary topological properties via finite topological spaces
}

First we must introduce notation for maps of finite topological spaces
we use. Two features are important for us:
\begin{itemize}
\item[ 1.] it reminds one that a finite topological space is a category
    (degenerate if you like)
\item[ 2.] it does not allow one to talk conveniently about non-identity
    endomorphisms of finite topological spaces. We hope this may help define 
    a decidable fragment of elementary topology because there is 
    a decidable fragment of diagram chasing without endomorphisms, see [GLZ].
\end{itemize}

A topological space comes with a {\em specialisation preorder} on its points: for
points $x,y \in X$,  $x \leq y$ iff $y \in cl x$ , or equivalently,  a category whose
objects are points of \ensuremath{X} and there is a unique morphism $x{\small\searrow}y$ iff $y \in cl x$.

For a finite topological space X, the specialisation preorder or equivalently
the category uniquely determines the space: a subset of \ensuremath{X} is closed iff it is
downward closed, or equivalently, there are no morphisms going outside the
subset.

The monotone maps (i.e. functors) are the continuous maps for this topology.

We denote a finite topological space by a list of the arrows (morphisms) in
the corresponding category; '$\leftrightarrow $' denotes an isomorphism and '$=$' denotes
the identity morphism.  An arrow between two such lists denotes
a continuous map (a functor) which sends each point to the correspondingly
labelled point, but possibly turning some morphisms into identity morphisms,
thus gluing some points. 

Thus, each point goes to "itself" and  $$
     \{a,b\}\longrightarrow \{a{\small\searrow}b\}\longrightarrow \{a\leftrightarrow b\}\longrightarrow \{a=b\}
$$
denotes
$$
   (discrete\ space\ on\ two\ points)\longrightarrow (Sierpinski\ space)\longrightarrow (antidiscrete\ space)\longrightarrow (single\ point)
$$

In $A \longrightarrow  B$, each object and each morphism in $A$ necessarily appears in \ensuremath{B} as well; we avoid listing 
the same object or morphism twice. Thus 
both 
$$
\{a\} \longrightarrow  \{a,b\}\text{ and } \{a\} \longrightarrow  \{b\}
$$  denote the same map from a single point to the discrete space with two points.
Both 
 $$\{a{\small\swarrow}U{\small\searrow}x{\small\swarrow}V{\small\searrow}b\}\longrightarrow \{a{\small\swarrow}U=x=V{\small\searrow}b\}\text{ and }\{a{\small\swarrow}U{\small\searrow}x{\small\swarrow}V{\small\searrow}b\}\longrightarrow \{U=x=V\}$$
denote the morphism gluing points $U,x,V$.

In $\{a{\small\searrow}b\}$, the point $a$ is open and point \ensuremath{b} is closed.

Let $$
  C_{<n} := \{ \ensuremath{f} : \ensuremath{f} \in \ensuremath{C},\text{ both the domain and range of }f \text{ are finite of size
                  less than }n \}.
$$

\begin{claim} 
  The following is a list of properties defined using the lifting
  property starting from a single morphism between spaces of at most two points.

  In the category of topological spaces, it holds:
\begin{itemize}
\item[]$  ((\{a\}\longrightarrow \{a{\small\searrow}b\})^r_{<5})^{lr}$ is almost the class of proper maps, namely
                            a map of T4 spaces is in the class iff it is proper

\item[]$   (\{b\}\longrightarrow \{a{\small\searrow}b\})^l  $        is the class of maps with dense image
\item[]$   (\{b\}\longrightarrow \{a{\small\searrow}b\})^{lr} $        is the class of maps of closed inclusions $A \subset  X$, $A$ is closed

\item[]$   (\emptyset\longrightarrow \{*\})^r=(\{0\}\longrightarrow \{0\leftrightarrow 1\})^l $  is the class of surjections
\item[]$   (\emptyset\longrightarrow \{*\})^{rl}      $      is the class of maps of form $A\longrightarrow A\cup D$, $D$ is discrete
\item[]$   (\emptyset\longrightarrow \{*\})^{rll} = (\{a\}\longrightarrow \{a,b\})^l$ is the class of maps $A\longrightarrow B$ such that
                             each open closed non-empty subset  
                             of \ensuremath{B} intersects $Im A$.
       
\item[]$   (\emptyset\longrightarrow \{*\})^l$ is the class of maps $A\longrightarrow B$ such that $A=B=\emptyset$ or $A\neq \emptyset$, $B$ arbitrary
\item[]$   (\emptyset\longrightarrow \{*\})^{ll}$ is the class of maps $A\longrightarrow B$ such that either $A=\emptyset$ or the map is an isomorphism
\item[]$   (\emptyset\longrightarrow \{*\})^{lll}$ is the class of maps $A\longrightarrow B$ which split

\item[]$   (\emptyset\longrightarrow \{*\})^{rr}$ is the class of subsets, i.e. injective maps $A\hookrightarrow B$ where the topology on A is induced from B.

\item[]$   (\{a\leftrightarrow b\}\longrightarrow \{a=b\})^l$      is the class of injections
\item[]$   (\{a{\small\searrow}b\}\longrightarrow \{a=b\})^l$       is the class of maps $f:X\longrightarrow Y$
                              such that the topology on $X$ is induced from $Y$
\item[]$   (\{a,b\}\longrightarrow \{a=b\})^l $       describes being connected, and is the class of maps $f:X\longrightarrow Y$ such that
                              $f(U) \cap f(V)=\emptyset$ for each two open closed subsets $U\neq V$ of $X$; 
                              if both $X$ and $Y$ are unions of open closed connected subsets, 
                              this means that the map  $\pi_0(X) \hookrightarrow  \pi_0(Y)$ is injective
\item[]$   (\{a\leftrightarrow b\}\longrightarrow \{a=b\})^r$      fibres are T0 spaces
\item[]$   (\{a{\small\searrow}b\}\longrightarrow \{a=b\})^r$       fibres are T1 spaces
\item[]$   (\{a,b\}\longrightarrow \{a=b\})^r $       is the class of injections

\item[]$   (\{a\}\longrightarrow \{a\leftrightarrow b\})^l $     is the class of surjections


\item[]  $ (\{a\}\longrightarrow \{a\leftrightarrow b\})^r $     is the class of surjections 
\item[]  $ (\{b\}\longrightarrow \{a{\small\searrow}b\})^l $     something T1-related but not particularly nice
\item[]  $ (\{a\}\longrightarrow \{a{\small\searrow}b\})^l$      something T0-related
\item[]  $ (\{a\}\longrightarrow \{a,b\})^l$        is the class of maps $f:X\longrightarrow Y$
                            such that either $X$ is empty or $f$ is surjective
\end{itemize}
\end{claim}
{\bf
Proof.} All items are trivial to verify, with the possible exception of the first
item.  [Bourbaki, General Topology, I\S10.2, Thm.1(d), p.101], quoted in Appendix B, gives a characterisation
of proper maps by a lifting property with respect to maps associated to
ultrafilters. Using this it is easy to check that each map in $(\{a\}\longrightarrow \{a{\small\searrow}b\})^r_{<5}$ 
being closed, hence proper, implies that each map in $((\{a\}\longrightarrow \{a{\small\searrow}b\})^r_{<5})^{lr}$ is proper.
A theorem of [Taimanov], cf.~[Engelking, 3.2.1,p.136], quoted in Appendix B, 
states that for a compact Hausdorff space $K$,
a  Hausdorff space $K$ is compact iff the map $K\longrightarrow \{*\}$ is in  $C_T^{lr}$ where
$$ C_T:=    \{\, \{a\leftrightarrow b\}\longrightarrow \{a=b\},\, \{a{\small\searrow}b\}\longrightarrow \{a=b\},\,
     \{b\}\longrightarrow \{a{\small\searrow}b\},\,\{a{\small\swarrow}o{\small\searrow}b\}\longrightarrow \{a=o=b\}\,\,\}$$
It is easy to check that all the maps listed in the formula above are closed, hence proper, 
and therefore 
$$ C_T^{lr} 
\subseteq 
((\{a\}\longrightarrow \{a{\small\searrow}b\})^r_{<5})^{lr}$$  
Finally, note that the proof of Taimanov theorem generalises to
give that a proper map between normal Hausdorff (T4) spaces is in the larger class.
 QED.

\begin{conjecture}
      In the category of topological spaces,
      $$((\{a\}\longrightarrow \{a{\small\searrow}b\})^r_{<5})^{lr}$$ is the class of proper maps.
\end{conjecture}

\begin{remark}
    It is easy to see that $((\{a\}\longrightarrow \{a{\small\searrow}b\})^r_{<m})^{lr} \subset   ((\{a\}\longrightarrow \{a{\small\searrow}b\})^r_{<n})^{lr}$
    for any $m<n$. However, I do not know whether there is $n>m>3$ such that the inclusion is strict.
    An example using cofinite topology (suggested by Sergei Kryzhevich) 
    shows that $C_T^{lr}$
   does not define the class of compact spaces: indeed, consider infinite sets $A\subset B$, $\omega\leq \card A<\card B$, equipped
   with cofinite topology (i.e. a subset is closed iff it is finite). 
   Then $A\subseteq \ensuremath{B} \in C_T^l$ yet $A \subseteq \ensuremath{B} \,\rightthreetimes\,  A\longrightarrow \{*\}$ fails: 
    for a map $f:B\longrightarrow A$ the preimage of some (necessarily closed) point is infinite  as $\card \ensuremath{B} > \card A$, hence not closed,
   and the map is not continious. Hence, $A\longrightarrow \{*\} \,\notin\, C_T$ yet $A$ is compact (non-Hausdorff).
   This example could probably be generalised to show that
    that  $((\{a\}\longrightarrow \{a{\small\searrow}b\})^r_{<4})^{lr} \subsetneq   ((\{a\}\longrightarrow \{a{\small\searrow}b\})^r_{<5})^{lr} $.
\end{remark}

\begin{problem}
\begin{itemize}
\item[(a)]Calculate $$\!\!\!\!\!\!\!\!\!\!\!\! ((\{b\}\longrightarrow \{a{\small\searrow}b\})^r_{<5})^{lr},\  ((\{b\}\longrightarrow \{a{\small\searrow}b\})^{lrr},\ \text{and}\ 
   (\{a{\small\swarrow}U{\small\searrow}x{\small\swarrow}V{\small\searrow}b\}\longrightarrow \{a{\small\swarrow}U=x=V{\small\searrow}b\})^{lr} $$
   Could either be viewed as a "definition" of the real line?
\item[(b)]Characterise the interval $[0,1]$, a circle $\mathbb S^1$ and, more generally,  spheres $\mathbb S^n$ using their topological characterisations 
provided by the Kline sphere charterisation theorem and its analogues. An example of such a characterisation is
that a topological space $X$ is homomorphic to the circle $\mathbb S^1$ iff 
$X$ is a connected Hausdorff metrizable space such that $X\setminus \{x,y\}$ is not connected for any two points $x\neq y\in X$
([Hocking,Young, Topology, Thm.2-28,p.55]); another example is that a topological space $X$ is homomorphic to 
the closed interval $[0,1]$ iff 
$X$ is a connected Hausdorff metrizable space such that $X\setminus \{x\}$ is not connected for exactly two points $x\neq y\in X$
([Hocking,Young, Topology, Thm.2-27,p.54]).
\end{itemize}
\end{problem}
\begin{remark}
    Urysohn lemma and Tietze extension theorem relate lifting properties involving $\RR$ 
   and those involving opens maps of finite topological spaces, 
   and this is why we hope the question above might have something to
    do with the real line. Let us give some more details.

    Note a map \ensuremath{f} of finite spaces is open iff \ensuremath{f} is in $(\{b\}\longrightarrow \{a{\small\searrow}b\})^r$,
    and that $\{a{\small\swarrow}U{\small\searrow}x{\small\swarrow}V{\small\searrow}b\}\longrightarrow \{a{\small\swarrow}U=x=V{\small\searrow}b\}$ is an open map.

   A space \ensuremath{X} is normal (T4) iff $\emptyset\longrightarrow \ensuremath{X} \,\rightthreetimes\,  \{a{\small\swarrow}U{\small\searrow}x{\small\swarrow}V{\small\searrow}b\}\longrightarrow \{a{\small\swarrow}U=x=V{\small\searrow}b\}$, hence
   the Uryhson lemma can be stated as follows:
   $$
\emptyset\longrightarrow \ensuremath{X} \,\rightthreetimes\,  \{a{\small\swarrow}U{\small\searrow}x{\small\swarrow}V{\small\searrow}b\}\longrightarrow \{a{\small\swarrow}U=x=V{\small\searrow}b\}$$ iff 
$$\emptyset\longrightarrow \ensuremath{X} \,\rightthreetimes\,  \{0'\} \cup  [0,1] \cup  \{1'\} \longrightarrow  \{0=0'{\small\searrow}x{\small\swarrow}1=1'\}$$
   where points $0',0$ and $1,1'$ are topologically indistinguishable in
    $\{0'\} \cup  [0,1] \cup  \{1'\}$,
    the closed interval $[0,1]$ goes to $x$, and $0,0'$ map to point $0=0'$,  and $1,1'$
    map to point $1=1'$.

    Tietze extension theorem states that for a normal space \ensuremath{X} and a closed subset A of X,
         $A\longrightarrow \ensuremath{X} \,\rightthreetimes\,  \RR\longrightarrow \{*\} $, i.e. in notation
         $\RR\longrightarrow \{*\}$ is in  $$\left( (\{b\}\longrightarrow \{a{\small\searrow}b\})^{lr}\cap \{ A\longrightarrow \ensuremath{X} : \emptyset\longrightarrow \ensuremath{X} \,\rightthreetimes\,  \{a{\small\swarrow}U{\small\searrow}x{\small\swarrow}V{\small\searrow}b\}\longrightarrow \{a{\small\swarrow}U=x=V{\small\searrow}b\} \right)^r$$

   Note that $\emptyset\longrightarrow \ensuremath{X} \,\rightthreetimes\,  [0,1]\longrightarrow \{0{\small\searrow}x{\small\swarrow}1\}$ is the stronger property \ensuremath{X} is
   perfectly normal.  A normal space \ensuremath{X} is perfectly normal provided
   each closed subset of \ensuremath{X} is the intersection of a countably many open subsets.

   For example, is $\{0'\} \cup  [0,1] \cup  \{1'\} \longrightarrow  \{0=0'{\small\searrow}x{\small\swarrow}1=1'\}$ in
     $(\{a{\small\swarrow}U{\small\searrow}x{\small\swarrow}V{\small\searrow}b\}\longrightarrow \{a{\small\swarrow}U=x=V{\small\searrow}b\}^{lr} \subset 
 ((\{b\}\longrightarrow \{a{\small\searrow}b\})^r_{<5})^{lr} $?
\end{remark}

\begin{remark}
Is there a model category or a factorisation system of interest
    associated with any of these lifting properties, for example compactness/properness?
\end{remark}

Many of the separation axioms can  be expressed as lifting properties with respect to maps
involving up to \ensuremath{4} points and the real line, see [Appendix A].

\section{
Elementary topology as diagram chasing computations with finite categories
}

Early works talk of topology in terms of {\em neighbourhood} systems $U_x$
where $U_x$ varies though 
 {\em open neighbourhoods of points} of a topological space; this is how
the notion of topology was defined by Hausdorff. In the notation of arrows, 
{\em a neighbourhood system $U_x$, $x\in X$} would correspond to a system of arrows
$$
\{x\}\longrightarrow X\xrightarrow{\,\,\, U\,\,\,}\{x{\small\searrow}x'\}$$
and Hausdorff's axioms (A),(B),(C) (see Appendix~B) would correspond to 
diagram chasing rules. 

Here we show the axioms of topology stated in the more modern
language of open subsets 
can be seen as diagram chasing rules for manipulating diagrams
involving notation such as
$$   \{x\}\longrightarrow X,\, X\longrightarrow \{x{\small\searrow}y\},\, X\longrightarrow \{x\leftrightarrow y\}   $$
in the following straightforward way; cf. [Gavrilovich, Elementary Topology,\S.2.1] for more
details.

As is standard in category theory, identify a point $x$ of a topological space $X$
with the arrow $\{x\}\longrightarrow X$, a subset $Z$ of $X$ with the arrow $X\longrightarrow \{z\leftrightarrow z'\}$,
and an open subset $U$ of $X$ with the arrow $X\longrightarrow \{u{\small\searrow}u'\}$.
With these identifications, the Hausdorff axioms of a topological space become
rules for manipulating such arrows, as follows.

{\em Both the empty set and the whole of \ensuremath{X} are open} says that the compositions
$$ X\longrightarrow \{c\}\longrightarrow \{o{\small\searrow}c\}\text{ and }X\longrightarrow \{o\}\longrightarrow \{o{\small\searrow}c\}  $$
behave as expected (the preimage of \{o\} is empty under the first map,
and is the whole of \ensuremath{X} under the second map).

{\em The intersection of two open subsets is open} means the arrow
$$   X\longrightarrow \{o{\small\searrow}c\}\times\{o'{\small\searrow}c'\} $$
 behaves as expected (the ``two open subsets'' are the preimages of points $o\in\{o{\small\searrow}c\}$ and $o'\in\{o'{\small\searrow}c'\}$;
``the intersection'' is the preimage of $(o,o')$ in
$\{o{\small\searrow}c\}\times\{o'{\small\searrow}c'\}$ ).

Finally, {\em a subset $U$ of $X$ is open iff each point $u$ of $U$ has an open
neighbourhood inside of $U$}
 corresponds to the following diagram chasing rule:


for each arrow $X\xra[\xi_U]{}  \{U\longleftrightarrow \bar U \}$ it holds\\
 $ \xymatrix{ { }  & \{U\ra\bar U\} \ar[d] \\
{X} \ar[r]_{\!\!\!\!\!\!\!\!\!\xi_U} \ar@{-->}[ur] & { \{U\longleftrightarrow \bar U \} }}$
\ \ \ \ iff  for each $\{u\}\lra X$, \ \ \ \ \ \ \
 $ \xymatrix{ {\{u\}} \ar[r] \ar[d] & \{u \ra U \longleftrightarrow \bar U\} \ar[d] \\
{X} \ar[r]_{\!\!\!\!\!\!\!\!\!\!\!\!\!\!\!\!\!\!\xi_U} \ar@{-->}[ur] & { \{u=U\longleftrightarrow \bar U \} }}$

%
%
%
%
%
 {\em The preimage of an open set is open} corresponds to the composition $$
 X\longrightarrow Y\longrightarrow \{u{\small\searrow}u'\}\longrightarrow \{u\leftrightarrow u'\} .                                          $$

This observation suggests that some arguments in elementary topology may be
understood entirely in terms of diagram chasing, see [Gavrilovich, Elementary Topology] for some
examples. This reinterpretation may help clarify the nature of the axioms of a
topological space, in particular it offers a constructive approach, may clarify
to what extent set-theoretic language is necessary, and perhaps help to suggest
an approach to ''tame topology'' of Grothendieck, Does this lead to tame
topology of Grothendieck, i.e. a foundation of topology "without false
problems" and "wild phenomena" "at the very beginning" ?

Let us state two problems; we hope they help to clarify the notion
of an ergosystem and that of a topological space.

\begin{problem}
\label{problem:1}
Write a short program which extracts diagram chasing derivations
from texts on elementary topology, in the spirit of the ideology of
ergosystems/ergostructures.
\end{problem}

\begin{problem}\label{problem:2} 
Develop elementary topology in terms of  finite categories (viewed
as finite topological spaces) and labelled commutative diagrams, with an
emphasis on labels (properties) of morphisms defined by the Quillen lifting
property.  Does this lead to tame topology of Grothendieck, i.e. a foundation
of topology "without false problems" and "wild phenomena" "at the very
beginning" ?
\end{problem}

\begin{problem}
(Ganesalingam, Gowers) wrote an automatic theorem prover trying to
make it "thinking in a human way". In a couple of their examples their
generated proofs amount to diagram chasing, e.g.  Claim (v) shows the generated
proof of the claim that a closed subspace of a complete metric space is
necessarily complete translates to two applications of a diagram chasing rule
corresponding to the lifting property, see [Gavrilovich, Slides] for details.
Arguably, both examples on the first page of (ibid.) also correspond to diagram
chasing.  In their approach, can this be seen as evidence that humans are
really thinking by diagram chasing and the lifting property in particular? Note
that Claim (v) involves examples typically shown to students to clarify the
concepts of a metric space being complete or closed.  Can one base a similar
theorem prover on our observations, particularly in elementary topology?
\end{problem}

Let us comment on an approach to Problem 1.

We observed that there is a simple rule which leads to several notions in
topology interesting enough to be introduced in an elementary course. Can this
rule be extended to a very short program which learns elementary topology?

We suggest the following naive approach is worth thinking about.

The program maintains a collection of directed labelled graphs and certain
distinguished subgraphs. Directed graphs represent parts of a category;
distinguished subgraphs represent commutative diagrams.  Labels represent
properties of morphisms.  Further, the program maintains a collection of rules
to manipulate these data, e.g. to add or remove arrows and labels.

The program interacts with a flow of signals, say the text
of [Bourbaki, General Topology, Ch.1], and seeks correlations
between the diagram chasing rules and the flow of signals.
It finds a "correlation" iff certain strings occur nearby
in the signal flow iff they occur nearby in a diagram chasing rule.
To find "what's interesting",  by brute force it searches
for a valid derivation which exhibits such correlations.  To guide the search
and exhibit missing correlations in a derivation under consideration, it may
ask questions: are these two strings related?  Once it finds such a derivation,
the program "uses it for building its own structure".  Labels correspond to
properties of morphisms.  Labels defined by the lifting property play an
important role, often used to exclude counterexamples making a diagram chasing
argument fail.  In [DeMorgan] we analysed the text of the definitions of
surjective and injective maps showing what such a correlation may look like in
a "baby" case.

A related but easier task is to write a theorem prover doing diagram chasing 
in a model category. The axioms of a (closed or not) model category 
as stated in [Quillen,I.1.1] can be interpreted as rules to manipulate 
labelled commutative diagrams in a labelled category. It appears
straightforward how to formulate a logic (proof system) based on these rules
which would allow to express statements like: Given a labelled commutative diagram, 
(it is permissible to) add this or that arrow or label.
Moreover, it appears not hard to write a theorem prover for this logic   
doing brute force guided search. What is not clear whether this logic is complete
in any sense or whether there are non-trivial inferences of this form to prove. 

Writing such a theorem prover is particularly easy when the underlying category
of the model category is a partial order [GavrilovichHasson] and [BaysQuilder] wrote 
some code for doing diagram chasing in such a category. However, the latter
problem is particularly severe as well.

The two problems  are related; we hope they help to clarify the notion
of an ergosystem and that of a topological space.

The following is a more concrete question towards Problem~\ref{problem:2}.

\begin{problem}
   (a) Prove that a compact Hausdorff space is normal by diagram chasing;
       does it require additional axioms?
       Note that we know how to express the statement entirely in terms of the
       lifting property and finite topological spaces of small size.

   (b) Formalise the argument in [Fox, 1945] which implies the category of
       topological spaces is not Cartesian closed. Namely, Theorem \ensuremath{3} [ibid.] proves that
          if \ensuremath{X} is separable metrizable space, \ensuremath{R} is the real line, then
           \ensuremath{X} is locally compact
            iff
           there is a topology on $X^\RR$ such that for any space $T$,
           a function $h:X\times T\longrightarrow \RR$ is continuous iff
           the corresponding function $h^*:T\longrightarrow X^\RR$ is continuous
           (where $h(x,t)=h^*(t)(x))$
     Note that here we do not know how to express the statement.
\end{problem}

\begin{remark}\label{simp:obj:unifrom}
 In [Gavrilovich, Tame Topology,\S6] and [Gavrilovich, Elementary Topology,\S2.9] we note that
  uniform spaces may be viewed as simplicial objects in
  the category of topological spaces in the following way.
  We wish to emphasise that this observation can be easily "read off"
  from [Bourbaki, General Topology,  II\S1.1.1] if one is inclined
  to translate everything into diagram chasing.

  Let \ensuremath{X} be a set.  Let
  $$X\Longleftrightarrow X\times \ensuremath{X} \Longleftrightarrow X\times X\times \ensuremath{X} ...$$  
be the "trivial" simplicial set where degeneracy and face maps are
  coordinate projections and diagonal embeddings.
  A uniform structure on a set $X$ is a filter, hence a topology, on the set \ensuremath{X\times X}
  satisfying certain properties; equip \ensuremath{X\times X} with this topology.
  Put the antidiscrete topology on X; put on \ensuremath{X\times X\times X} the topology which is
  the pullback of the topologies on \ensuremath{X\times X} and \ensuremath{X} along the projection maps, and
  similarly for \ensuremath{X\times X\times X\times X} etc. A straightforward verification shows this is well-defined
  whenever the topology on \ensuremath{X\times X} corresponds to a uniform structure on X.
\end{remark}

\begin{remark}
In [Gavrilovich, Tame Topology, \S5.4] we observe that a number of consequences
   of compactness can be expressed as a change of order of quantifiers in a
   formula,
   i.e. are of form 
   $\forall \exists \phi(...)\implies\exists\forall \phi(...)$
   namely that a real-valued function on a compact is necessarily bounded,
   that a Hausdorff compact is necessarily normal,
   that the image in \ensuremath{X} of a closed subset in $X\times K$ is necessarily closed,
   Lebesgue's number Lemma, and paracompactness.

   Such formulae correspond to inference rules of a special form,
   and we feel a special syntax should be introduced to state
   these rules.

   For example, consider the statement that
   "a real-valued function on a compact domain is necessarily bounded".
   As a first order formula, it is expressed as
$$
   \forall \ensuremath{x} \in \ensuremath{K} \exists M ( f(x) \leq  M ) \Longrightarrow  \exists M \forall \ensuremath{x} \in \ensuremath{K} ( f(x) \leq  M )
$$
   Another way to express it is:
$$
   \exists M:K\longrightarrow \RR \forall \ensuremath{x} \in \ensuremath{K} ( f(x) \leq  M(x) ) \Longrightarrow  \exists M \in \RR \forall \ensuremath{x} \in \ensuremath{K} ( f(x) \leq  M )
$$
   Note that all that happened here is that a function $M:K\longrightarrow \RR$,
   become a constant $M \in \RR$, or rather
   expression "M(x)" of type $K\longrightarrow M$ which used (depended upon) variable "x"
   become expression "M" which does not use (depend upon)  variable "x".

   We feel there should be a special syntax which would allow
   to express above as an inference rule {\em removing dependency of "M(x)" on "x"},
   and this syntax should be used to express consequences of compactness
   in a diagram chasing derivation system for elementary topology.

   To summarise, we think that compactness should be formulated
   with help of inference rules for expressly manipulating which variables are 'new', 
   in what order they 'were' introduced, 
   and what variables terms depend on, e.g. rules replacing a term t(x,y) by term t(x).

   Something like the following:

\begin{verbatim}
   ... f(x) =< M(x) ...
   --------------------
   ...  f(x) =< M   ...
\end{verbatim}
\end{remark}

\section{
  Ergo-Structures/Ergo-Systems Conjecture of Misha Gromov.
}

We conclude with a section which aims to explain our motivation and
hence is speculative and perhaps somewhat inappropriate in what is mostly
a mathematical text.

Misha Gromov [Memorandum Ergo] conjectures
there are particular mathematical, essentially combinatorial, structures,
called {\em ergostructures} or {\em ergosystems},
which help to understand complex biological behaviour including learning
and create  mathematically interesting models of these processes.

We hope our observations may eventually help to uncover
an essentially combinatorial reasoning behind elementary topology,
and thereby suggest an example of an ergostructure.

An ergosystem/ergostructure may be thought of as an "engine" producing
(structurally or mathematically) {\em interesting} behaviour
which is then misappropriated into a {\em useful} behaviour by a biological
system.

"Behaviour"  is thought of as an interaction with a flow of signals.
By itself, such an "engine" produces {\em interesting} behaviour,
with little or no concern for any later use;
it "interacts with an incoming flow of signals;
it recognizes and selects what is {\em interesting} for itself
in such a flow and uses it for building its own structure"
[Gromov, Memorandum Ergo].

An analogy might help. Consider a complex mechanical contraption powered by an
engine, such as a loom. By itself, there is nothing directly useful done by its
engine; indeed, the very same engine may be used by different mechanisms for
all sorts of useful and useless tasks.  To understand the workings of a
mechanism, sometimes you had better forget its purpose and ask what is the
engine, how is it powered,  and what keeps the engine in good condition.
Understanding  the loom's engine (only) in terms of how it helps to weave is
misleading.

Thus, the concept of an ergostructure/ergosystems suggests a different kind of
questions we should ask about biological systems and learning.

A further suggestion is that these "engines" might be rather universal, i.e.
able to behave interestingly interacting with a diverse range of signal flows.
At the very beginning an ergosystem/structure is a "crisp" mathematically
interesting structure of size small enough to be stumbled upon by evolution; as
it grows, it becomes "fuzzy" and specialised to a particular kind of  flow of
signals.

However, we want to draw attention to the following specific suggestion:
\begin{quote}
   ""The category/functor modulated structures can not be directly used by ergosystems,
     e.g. because the morphisms sets between even moderate objects are
     usually unlistable.
     But the ideas of the category theory show that there are certain (often
     non-obviuos) rules for generating proper concepts. (You ergobrain would not
     function if it had followed the motto: "in my theory I use whichever definitions
     I like".) The category theory provides a (rough at this stage) hint on a possible
     nature of such rules.
   [Gromov, Ergobrain]
\end{quote}
Our observations give an example of a simple rule which can be used "to
generate proper concepts", particularly in elementary topology. We hope that
our observations can make the hint less rough, particularly if one properly
develops elementary topology in terms of diagram chasing, with an emphasis on
the lifting property.

{\bf Problem 1.} Write a short program which extracts diagram chasing derivations
from texts on elementary topology, in the spirit of the ideology of
ergosystems/ergostructures.  That is, it considers a flow of signals
interesting if it correlates with diagram chasing rules.

In the previous section we give some suggestions, albeit naive, what such a program
might look like.

{\bf Problem 2.} Develop elementary topology in terms of labelled commutative
diagrams involving finite categories (viewed as finite topological spaces),
with an emphasis on labels (properties) of morphisms defined by the Quillen
lifting property.

 Does this lead to the tame topology of Grothendieck, i.e. a foundation of topology
"without false problems" and "wild phenomena" "at the very beginning" ?

In the previous section we give some suggestions, albeit naive, what such
a program might look like and how to express elementary topology
in terms of labelled diagram chasing.

Acknowledgements. To be written.

 This work is a continuation of [DMG]; early history is given there. I thank
M.Bays, K.Pimenov, V.Sosnilo, S.Synchuk and P.Zusmanovich for discussions and proofreading; I
thank L.Beklemishev, N.Durov, S.V.Ivanov, S.Podkorytov, A.L.Smirnov for discussions. I also
thank several students for encouraging and helpful discussions.  Chebyshev
laboratory, St.Petersburg State University, provided a coffee machine and an
excellent company around it to chat about mathematics.  Special thanks are to
Martin Bays for many corrections and helpful discussions. Several observations
in this paper are due to Martin Bays. I thank S.V.Ivanov for several
encouraging and useful discussions; in particular, he suggested to look at the
Lebesque's number lemma and the Arzela-Ascoli theorem. A discussion with Sergei
Kryzhevich motivated the group theory examples.

   Much of this paper was done in St.Petersburg; it wouldn't have been possible
without support of family and friends who created an excellent social
environment and who occasionally accepted an invitation for a walk or a coffee
or extended an invitation; alas, I made such a poor use of it all.

  This note is elementary, and it was embarrassing and boring, and
embarrassingly boring, to think or talk about matters so trivial, but luckily
I had no obligations for a time.

\newpage

\section{
 Appendix A. Separation axioms as lifting properties (from Wikipedia)
}

The separation axioms are lifting properties with respect to maps involving up to \ensuremath{4} points
and the real line. What follows below is the text of the Wikipedia page on the separation axioms
where we added lifting properties formulae expressing what is said there in words.

%

Let \ensuremath{X} be a topological space. Then two points \ensuremath{x} and \ensuremath{y} in \ensuremath{X} are {\em topologically distinguishable}
iff the map $\{x\leftrightarrow y\} \longrightarrow  X$ is not continuous, i.e. 
iff 
at least one of them has an open neighbourhood which is not a neighbourhood of the other.

Two points \ensuremath{x} and \ensuremath{y} are {\em separated} iff neither $ \{x{\small\searrow}y\} \longrightarrow  X$ nor $\{x{\small\searrow}y\} \longrightarrow  X$ is continuous, 
i.e~each of them has a neighbourhood that is not a neighbourhood of the other; 
in other words, neither belongs to the other's closure, $x \notin cl\, x$ and $y \notin cl\, x$. 
More generally, two subsets A and \ensuremath{B} of \ensuremath{X} are {\em separated} iff each is disjoint from the other's closure, 
i.e.~$ A\cap cl \ensuremath{B} = B\cap cl A = \emptyset $. 
(The closures themselves do not have to be disjoint.) In other words, the map
$ i_{AB} : \ensuremath{X} \longrightarrow  \{A \leftrightarrow  \ensuremath{x} \leftrightarrow  B\}$
sending the subset $A$ to the point $A$, the subset $B$ to the point $B$, and the rest to the point $x$,  
factors both as 
 $$ 
X \longrightarrow   \{A \leftrightarrow  U_A {\small\searrow} \ensuremath{x} \leftrightarrow  B\} \longrightarrow  \{A=U_A \leftrightarrow  \ensuremath{x} \leftrightarrow  B\} $$
and 
 $$ 
X \longrightarrow   \{A \leftrightarrow  \ensuremath{x} {\small\swarrow} U_B \leftrightarrow  B\} \longrightarrow  \{A \leftrightarrow  \ensuremath{x} \leftrightarrow  U_B=B\} $$ 
here the preimage of $x,B$, resp. $x,A$ is a closed subset containing $B$, resp. $A$, and disjoint from $A$, resp. $B$.
All of the remaining conditions for separation of sets may also be applied to points (or to a point and a set) 
by using singleton sets. Points \ensuremath{x} and \ensuremath{y} will be considered separated, by neighbourhoods, 
by closed neighbourhoods, by a continuous function, precisely by a function, 
iff their singleton sets $\{x\}$ and $\{y\}$ are separated according to the corresponding criterion.

Subsets A and \ensuremath{B} are {\em separated by neighbourhoods} iff
A and \ensuremath{B} have disjoint neighbourhoods, i.e. 
iff  
$ i_{AB} : \ensuremath{X} \longrightarrow  \{A \leftrightarrow  \ensuremath{x} \leftrightarrow  B\}$  factors as $$ 
X \longrightarrow   \{A \leftrightarrow  U_A {\small\searrow} \ensuremath{x} {\small\swarrow} U_B \leftrightarrow  B\} \longrightarrow  \{A=U_A \leftrightarrow  \ensuremath{x} \leftrightarrow  U_B=B\} $$
here the disjoint neighbourhoods of A and \ensuremath{B} are the preimages of open subsets ${A,U_A}$ and ${U_B,B}$ of 
$ \{A \leftrightarrow  U_A {\small\searrow} \ensuremath{x} {\small\swarrow} U_B \leftrightarrow  B\}$,  resp. 
They are {\em separated by closed neighbourhoods} 
iff they have disjoint closed neighbourhoods, i.e.
$i_{AB}$ factors as 
$$
X \longrightarrow   \{A \leftrightarrow  U_A {\small\searrow} U'_A {\small\swarrow} \ensuremath{x} {\small\searrow} U'_B {\small\swarrow} U_B \leftrightarrow  B\} \longrightarrow  \{A\leftrightarrow U_A=U'_A= \ensuremath{x} = U'_B=U_B\leftrightarrow B\} 
.$$
They are {\em separated by a continuous function} iff 
there exists a continuous function \ensuremath{f} from the space \ensuremath{X} to the real line $\RR$ such that $f(A)=0$  and $f(B)=1$,
i.e.
the map $i_{AB}$ factors as
$$
X \longrightarrow   \{0'\} \union [0,1] \union \{1'\}  \longrightarrow  \{A \leftrightarrow  \ensuremath{x} \leftrightarrow  B\}  
$$
where points $0',0$ and $1,1'$ are topologically indistinguishable, 
and $0'$ maps to $A$, and $1'$ maps to $B$, and $[0,1]$ maps to $x$. 
Finally, they are {\em precisely separated by a continuous function}
iff there exists a continuous function \ensuremath{f} from \ensuremath{X} to $\RR$ such that the preimage $f^{ - 1}(\{0\})= A$ and $f^{ - 1}(\{1\})=B$.
i.e.~iff $i_{AB}$ factors as $$ 
X \longrightarrow   [0,1]  \longrightarrow  \{A \leftrightarrow  \ensuremath{x} \leftrightarrow  B\} 
$$ where $0$ goes to point $A$ and $1$ goes to point $B$.

These conditions are given in order of increasing strength: 
Any two topologically distinguishable points must be distinct, and any two separated points must be topologically distinguishable. Any two separated sets must be disjoint, any two sets separated by neighbourhoods must be separated, and so on.


The definitions below all use essentially the preliminary definitions above.

In all of the following definitions, \ensuremath{X} is again a topological space.
\begin{itemize}
\item[]    \ensuremath{X} is T0, or Kolmogorov, if any two distinct points in \ensuremath{X} are topologically
distinguishable. (It will be a common theme among the separation axioms to have
one version of an axiom that requires T0 and one version that doesn't.)
As a formula, this is expressed as
   $$ \{x\leftrightarrow y\} \longrightarrow  \{x=y\} \,\rightthreetimes\,  \ensuremath{X} \longrightarrow  \{*\}$$

\item[]    \ensuremath{X} is R0, or symmetric, if any two topologically distinguishable points in X
are separated, i.e.
         $$\{x{\small\searrow}y\} \longrightarrow  \{x\leftrightarrow y\} \,\rightthreetimes\,  \ensuremath{X} \longrightarrow  \{*\}$$

\item[]    \ensuremath{X} is T1, or accessible or Frechet, if any two distinct points in \ensuremath{X} are
separated, i.e.
$$     \{x{\small\searrow}y\} \longrightarrow  \{x=y\} \,\rightthreetimes\,  \ensuremath{X} \longrightarrow  \{*\}  $$
 Thus, \ensuremath{X} is T1 if and only if it is both T0 and R0. (Although you may
say such things as "T1 space", "Frechet topology", and "Suppose that the
topological space \ensuremath{X} is Frechet", avoid saying "Frechet space" in this
context, since there is another entirely different notion of Frechet space in
functional analysis.)

\item[]    \ensuremath{X} is R1, or preregular, if any two topologically distinguishable points in
X are separated by neighbourhoods. Every R1 space is also R0.

\item[]    \ensuremath{X} is weak Hausdorff, if the image of every continuous map from a compact
Hausdorff space into \ensuremath{X} is closed. All weak Hausdorff spaces are T1, and all
Hausdorff spaces are weak Hausdorff.

\item[]    \ensuremath{X} is Hausdorff, or T2 or separated, if any two distinct points in \ensuremath{X} are
separated by neighbourhoods, i.e. 
$$    \{x,y\} \hookrightarrow  \ensuremath{X} \,\rightthreetimes\,  \{x{\small\searrow}X{\small\swarrow}y\} \longrightarrow  \{x=X=y\} $$
Thus, \ensuremath{X} is Hausdorff if and only if it is both T0
and R1. Every Hausdorff space is also T1.

\item[]    \ensuremath{X} is $T2\frac12$, or Urysohn, if any two distinct points in \ensuremath{X} are separated by
closed neighbourhoods, i.e. 
$$
 \{x,y\} \hookrightarrow  \ensuremath{X} \,\rightthreetimes\,  \{x{\small\searrow}x'{\small\swarrow}X{\small\searrow}y'{\small\swarrow}y\} \longrightarrow  \{x=x'=X=y'=y\}
$$
Every T2$\frac12$ space is also Hausdorff.

\item[]    \ensuremath{X} is completely Hausdorff, or completely T2, if any two distinct points in
X are separated by a continuous function, i.e. 
$$     \{x,y\} \hookrightarrow  \ensuremath{X} \,\rightthreetimes\,   [0,1]\longrightarrow \{*\}
$$     where  $\{x,y\} \hookrightarrow  X$  runs through all injective maps from the discrete two
point space $\{x,y\}$.

Every completely Hausdorff space is
also T2$\frac 12 ½$.

\item[]    \ensuremath{X} is regular if, given any point \ensuremath{x} and closed subset $F$ in \ensuremath{X} such that \ensuremath{x} does
not belong to $F$, they are separated by neighbourhoods, i.e. 
$$    \{x\} \longrightarrow  \ensuremath{X} \,\rightthreetimes\,  \{x{\small\searrow}X{\small\swarrow}U{\small\searrow}F\} \longrightarrow  \{x=X=U{\small\searrow}F\}
$$
(In fact, in a regular
space, any such \ensuremath{x} and\ensuremath{F} will also be separated by closed neighbourhoods.) Every
regular space is also R1.

\item[]    \ensuremath{X} is regular Hausdorff, or T3, if it is both T0 and regular.[1] Every
regular Hausdorff space is also $T2\frac12$.

\item[]    \ensuremath{X} is completely regular if, given any point \ensuremath{x} and closed set $F$ in \ensuremath{X} such
that \ensuremath{x} does not belong to $F$, they are separated by a continuous function, i.e. 
$$
      \{x\} \longrightarrow  \ensuremath{X} \,\rightthreetimes\,  [0,1] \cup \{F\} \longrightarrow  \{x{\small\searrow}F\}
$$    where points $F$ and $1$ are topologically indistinguishable, $[0,1]$ goes to $x$,
and $F$ goes to $F$.

Every
completely regular space is also regular.

\item[]    \ensuremath{X} is Tychonoff, or T3$\frac12$, completely T3, or completely regular Hausdorff, if
it is both T0 and completely regular.[2] Every Tychonoff space is both regular
Hausdorff and completely Hausdorff.

\item[]    \ensuremath{X} is normal if any two disjoint closed subsets of \ensuremath{X} are separated by
neighbourhoods, i.e.
$$   \emptyset \longrightarrow \ensuremath{X} \,\rightthreetimes\,  \{x{\small\swarrow}x'{\small\searrow}X{\small\swarrow}y'{\small\searrow}y\} \longrightarrow  \{x{\small\swarrow}x'=X=y'{\small\searrow}y\}
$$
 In fact, by Urysohn lemma a space is normal if and only if any two disjoint
closed sets can be separated by a continuous function, i.e.
 $$   \emptyset \longrightarrow  \ensuremath{X} \,\rightthreetimes\,  \{0'\} \cup [0,1] \cup \{1'\} \longrightarrow  \{0=0'{\small\searrow}x{\small\swarrow}1=1'\} $$ 
where points
$0',0$ and $1,1'$ are topologically indistinguishable,
            $[0,1]$ goes to $x$, and both $0,0'$ map to point $0=0'$,  and both $1,1'$ map to point
$1=1'$.

\item[]    \ensuremath{X} is normal Hausdorff, or T4, if it is both T1 and normal. Every normal
Hausdorff space is both Tychonoff and normal regular.

\item[]    \ensuremath{X} is completely normal if any two separated sets $A$ and $B$ are separated by
neighbourhoods $U\supset A$ and $V\supset B$
      such that $U$ and $V$ do not intersect, i.e.????
 $$\emptyset \longrightarrow  \ensuremath{X} \,\rightthreetimes\,  \{X{\small\swarrow}A\leftrightarrow U{\small\searrow}U'{\small\swarrow}W{\small\searrow}V'{\small\swarrow}V\leftrightarrow B{\small\searrow}X\} \longrightarrow  \{U=U',V'=V\}$$
      Every completely normal space is also normal.

\item[]    \ensuremath{X} is perfectly normal if any two disjoint closed sets are precisely
separated by a continuous function, i.e. 
$$    \emptyset\longrightarrow \ensuremath{X} \,\rightthreetimes\,  [0,1]\longrightarrow \{0{\small\swarrow}X{\small\searrow}1\}
$$   where $(0,1)$ goes to the open point $X$, and $0$ goes to $0$, and $1$ goes to $1$.

Every perfectly normal space is also
completely normal.

\item[]  \ensuremath{X} is extremally disconnected if the closure of every open subset of $X$ is open, i.e.
  $$\emptyset\longrightarrow \ensuremath{X} \,\rightthreetimes\,  \{U{\small\searrow}Z',Z{\small\swarrow}V\}\longrightarrow \{U{\small\searrow}Z'=Z{\small\swarrow}V\}$$
  or equivalently
   $$\emptyset\longrightarrow \ensuremath{X} \,\rightthreetimes\,  \{U{\small\searrow}Z',Z{\small\swarrow}V\}\longrightarrow \{Z'=Z\}$$ 

\end{itemize}

\newpage

\section{Appendix B. Quotations from sources.}

For reader's convenience we quote here from the several sources we use.

[Bourbaki, General Topology, I\S10.2, Thm.1(d), p.101]:
\begin{quote} THEOREM I. Let $f: X\longrightarrow  Y$ be a continuous mapping. Then the following
four statements are equivalent:
\begin{itemize}
\item[a)] \ensuremath{f} is proper.
\item[b)] \ensuremath{f} is closed and $f^{-1} (y)$ is quasi-compact for each $y\in Y$.
\item[c)] If $\mathfrak F$ 
 is a filter on \ensuremath{X} and if $y \in Y$ is a cluster point of $f (\mathfrak F)$ then there
is a cluster point \ensuremath{x} of 
 such that $f (x) = y$.
\item[d)] If $\mathfrak U$ is an ultrafilter on \ensuremath{X} and if $y \in Y$ is a limit point of the ultrafilter
base $f (\mathfrak U)$, then there is a limit point \ensuremath{x} of $\mathfrak U$ such that $f (x) = y$.
\end{itemize}
\end{quote}
[Engelking, 3.2.1,p.136] (``compact'' below stands for ``compact Hausdorff''):
\begin{quote} 3.2.1. THEOREM. Let $A$ be a dense subspace of a topological space \ensuremath{X} and \ensuremath{f} a continuous
mapping of $A$ to a compact space $Y$. The mapping \ensuremath{f} has a continuous extension over \ensuremath{X} if and
only if for every pair $B_1,B_2$ of disjoint closed subsets of \ensuremath{Y} the inverse images $f^{-1}(B_1)$ 
and $f^{-1}(B_2)$ have disjoint closures in the space $X$.
\end{quote}
[Hausdorff, Set theory, \S40, p.259] (``$\varepsilon$'' stands  for ``$\in$'', and
``$ U_x V_x$'' stands for ``$ U_x \cap V_x$'') :

\begin{quote}
From
the theorems about open sets we derive the following properties of the
neighborhoods:
\begin{itemize}
\item[(A)]Every point $x$ has at least one neighborhood $U_x$; and $U_x$ always
contains $x$.
\item[(B)]For any two neighborhoods $U_x$ and $V_x$ of the same point, there
exists a third, $W_x \leq  U_x V_x$.
\item[(C)]Every point $y \varepsilon U_z$ has a neighborhood $U_y \leq U_x$.
\end{itemize}
It is now again possible to treat neighborhoods as unexplained 
concepts and to use them as our starting point, postulating Theorems (A),
(B), and (C) as neighborhood axioms$.^1$ Open sets $G$ are then defined
as sums of neighborhoods or as sets in which every point $x \varepsilon G$ has a
neighborhood $U_x \leq G$ (the null set included). Theorems (1), (2), and
(3) about open sets are then provable.

....

$\\{\\}^1$ Such a program was carried through in the first edition of this book. [Grund-
z\"ugeder Mengenlehre. (Leipzig, 1914; repr. New York, 1949),]
\end{quote}


\end{document}